\documentclass[11pt]{article}

\usepackage{amssymb,latexsym}
\usepackage{amsmath,amscd}
\usepackage{citesort}

\title{{ Semiclassical Geometry of Quantum Line Bundles and Morita 
Equivalence of Star Products}}

\author{{     
         Henrique Bursztyn\thanks{E-mail: henrique@math.berkeley.edu}
         \thanks{Current address:
         Mathematical Sciences Research Institute, 1000 Centennial Dr., Berkeley, CA 94720.}
        } \\[0.3cm]
        Department of Mathematics\\
        University of California, Berkeley\\
        CA, 94720 \\}

\date{January 2002}

\newcommand{\im} {{\mathrm i}}

\newcommand{\id}         {{\mathsf {id}}}

\newcommand{\Aut}        {{\mathrm {Aut}}}

\newcommand{\Hom}        {{\mathrm {Hom}}}

\newcommand{\End}        {{\mathrm {End}}}
\newcommand{\Def}        {{\mathrm {Def}}_{\mbox{\tiny{diff}}}}
\newcommand{\Defa}       {{\mathrm {Def}}}
\newcommand{\Der}        {{\mathrm {Der}}}
\newcommand{\Pic}        {{\mathrm {Pic}}}

\newcommand{\ch}         {\mbox{c}}
\newcommand{\cL}         {\mbox{c}_1}
\newcommand{\pch}        {\mbox{c}_1^{\pi}}
\newcommand{\Right}      {{R}}

\newcommand{\qR} {\boldsymbol{R}}

\newcommand{\qA} {\boldsymbol{\mathcal A}}

\newcommand{\qE} {\boldsymbol{\mathcal E}}
\newcommand{\qP} {\boldsymbol{P}}

\newcommand{\qT} {\boldsymbol{T}}
\newcommand{\BEA} {{\sideset{_{\scriptscriptstyle \mathcal B}}
                            {_{\scriptscriptstyle\mathcal A}}
                            {\operatorname{\mathcal E}}}}
\newcommand{\SER} {{\sideset{_{\scriptscriptstyle \mathcal{S}}}
                            {_{\scriptscriptstyle \mathcal{R}}}
                            {\operatorname{\mathcal E}}}}
\newcommand{\SE} {{\sideset{_{\scriptscriptstyle \mathcal{S}}}
                            {}
                            {\operatorname{\mathcal E}}}}
\newcommand{\ER} {{\sideset{}
                            {_{\scriptscriptstyle \mathcal{R}}}
                            {\operatorname{\mathcal E}}}}

\newcommand{\AEA} {{\sideset{_{\scriptscriptstyle \mathcal A}}
                            {_{\scriptscriptstyle\mathcal A}}
                            {\operatorname{{\mathcal E}}}}}

\newcommand{\RES} {{\sideset{_{\scriptscriptstyle \mathcal{R}}}
                            {_{\scriptscriptstyle \mathcal{S}}}
                            {\operatorname{\overline{\mathcal E}}}}}
\newcommand{\RER} {{\sideset{_{\scriptscriptstyle \mathcal{R}}}
                            {_{\scriptscriptstyle \mathcal{R}}}
                            {\operatorname{{\mathcal E}}}}}

\newcommand{\RE} {{\sideset{_{\scriptscriptstyle \mathcal{R}}}
                            {}
                            {\operatorname{\overline{\mathcal E}}}}}
\newcommand{\ES} {{\sideset{}
                            {_{\scriptscriptstyle \mathcal{S}}}
                            {\operatorname{\overline{\mathcal E}}}}}

\newcommand{\EA} {{\sideset{}{_{\scriptscriptstyle\mathcal A}}
                            {\operatorname{\mathcal E}}}}
\newcommand{\BE} {{\sideset{_{\scriptscriptstyle \mathcal B}}{}
                             {\operatorname{\mathcal E}}}}

\newcommand{\AAA} {{\sideset{_{\scriptscriptstyle \mathcal A}}
                            {_{\scriptscriptstyle \mathcal A}}
                            {\operatorname{{\mathcal A}}}}}

\newtheorem{lemma} {Lemma} [section]
\newtheorem{proposition} [lemma] {Proposition}
\newtheorem{theorem} [lemma] {Theorem}
\newtheorem{corollary} [lemma] {Corollary}
\newtheorem{definition}[lemma] {Definition}
\newtheorem{example}[lemma] {Example}
\newtheorem{remark}[lemma]{Remark}

\newenvironment{proof}{{\sc Proof:}}{{\hspace*{\fill} $\square$\\}}

\numberwithin{equation}{section}

\begin{document}
\maketitle
\begin{abstract}
In this paper we show how deformation quantization of line bundles over a
Poisson manifold $M$ produces a canonical action $\Phi$ of the Picard group
$\Pic(M)\cong H^2(M,\mathbb Z)$ on the moduli space of equivalence classes
of differential star products on $M$, $\Def(M)$.
The orbits of $\Phi$ characterize Morita equivalent star products on $M$. 
We describe the semiclassical limit of $\Phi$
in terms of the characteristic classes of star products by studying 
the semiclassical geometry of deformed line bundles.
\end{abstract}

\section{Introduction}\label{IntroSec}
The notion of ``representation equivalence'' of objects in a given category was
first made precise by Morita \cite{Morita58} in the context of unital rings.
Since then, the concept of Morita equivalence has been studied in
 many other settings,
including groupoids, operator algebras and Poisson manifolds (see 
\cite{Land2000} for a unified approach, with original references).

Connections between noncommutative geometry and $M$-theory \cite{ConnDoSch98} 
have shown that Morita equivalence is related to physical duality \cite{Sch98};
this motivated the study of the classification of quantum tori up to Morita
equivalence \cite{RiefSch}. One can think of (the algebra of functions on) quantum
tori as objects obtained from (the algebra of functions on) ordinary tori, equipped
with a constant Poisson structure, by means of \emph{strict deformation
quantization} \cite{Rief89}; the classification 
problem is to determine when constant Poisson
structures $\theta, \theta'$ on $T^n$ 
give rise to Morita equivalent quantum tori $T_{\theta} \sim T_{\theta'}$.

An analogous problem can be considered in the framework of formal deformation
quantization \cite{BFFLS78}
(see \cite{Ste98,Wei94} for surveys). 
In this approach to quantization, algebras of quantum observables
are defined by formal deformations \cite{Ger64}
 of algebras of classical observables. These
deformed algebra structures are called \emph{star products}, and they often
(but not always) arise as asymptotic expansions of strict quantizations (see
\cite[Sec.~4]{Gutt2000} and references therein). 
Unlike strict quantizations, star products have been shown to exist on
arbitrary Poisson manifolds 
\cite{Kon97b}. In this paper, we 
investigate the problem of which star products
define Morita equivalent deformed algebras. 

We show that 
deformation quantization of line bundles (see e.g. \cite{BuWa2000b})
over a Poisson manifold $(M,\pi)$
produces a canonical action
$$
\Phi: \Pic(M) \times \Def(M,\pi) \longrightarrow \Def(M,\pi)
$$
of the Picard group of $M$ 
on $\Def(M,\pi)$, the moduli space of equivalence classes of differential star
products on $(M,\pi)$. The action $\Phi$ characterizes Morita equivalent
star products on $(M,\pi)$:  $\star$ and $\star'$
 are Morita equivalent if and only if there exists a Poisson diffeomorphism
$\psi: M \longrightarrow M$ so that the equivalence
classes $[\psi^* (\star')]$ and $[\star]$ lie 
in the same $\Phi$-orbit.

We use the well-known descriptions of the set $\Def(M,\pi)$ (in terms
of Fedosov-Deligne's characteristic classes in the symplectic case
\cite{Deligne95,Fed94a,NT95b,BCG97,WX98}, and 
in terms of Kontsevich's classes of formal Poisson structures \cite{Kon97b}
in the case of arbitrary Poisson manifolds) to compute
the semiclassical limit of $\Phi$ explicitly.
This involves the study of the semiclassical limit of line bundle deformations
over Poisson manifolds. Just as the semiclassical limit of deformations
of the associative algebra
structure of $C^\infty(M)$ gives rise to Poisson structures on $M$
\cite[Sec.~19.4]{SilWein99}, the semiclassical
limit of deformed line bundles defines a geometric object on the underlying
classical line bundle: a contravariant connection \cite{Rui2000,Vais91}.
Contravariant connections are analogous to ordinary connections, but with
cotangent vectors playing the role of tangent ones. 
They define a characteristic
class on line bundles over Poisson manifolds, called the Poisson-Chern class.
We show explicitly how the semiclassical limit of $\Phi$ ``twists'' characteristic classes
of star products by Poisson-Chern classes. 

As a result, it follows that,
when $M$ is symplectic, with free $H^2(M,\mathbb Z)$, the action $\Phi$
is faithful, and one gets a parametrization of star products Morita
equivalent to a fixed one (up to isomorphism). The discussion also
provides an integrality obstruction for Morita equivalent star products
in general.
In the follow-up paper
\cite{BuWa2001}, this integrality condition is shown to be related to
Dirac's quantization condition for magnetic charges.
We remark that similar results hold for ``strongly'' Morita equivalent
hermitian star products in the sense of \cite{BuWa99a,BuWa2000} (see
\cite{BuWa2001}). 

The paper is organized as follows.
Section \ref{PrelimSec} recalls some of the necessary background:
Morita equivalence, deformations of associative algebras (star products),
deformations of finitely generated projective modules (and vector bundles),
and Poisson fibred algebra structures.
In Section \ref{ActionSec} we define the action $\Phi$ and show how it is related to
Morita equivalence.
In Section \ref{SemiClassSec} we study the semiclassical geometry of quantum line
bundles over Poisson manifolds, and show how contravariant connections arise
in this context.
In Section \ref{EquivSec} we describe the semiclassical limit of Morita equivalent
star products in terms of their characteristic classes.
We have included a summary of standard results on Poisson geometry in
Appendix \ref{AppA}.

\noindent {\bf Acknowledgments} I would like to thank Alan Weinstein and
Stefan Waldmann for many helpful discussions and suggestions. Some ideas
related to this paper can be found in \cite{JSW2001}; I thank the authors
for bringing their work to my attention and for discussions that led to 
corrections on a previous version of this note.

\section{Preliminaries}\label{PrelimSec}

\subsection{Morita equivalence and the Picard group}\label{MorPicSub}

For a unital ring $\mathcal R$, let $_{\scriptscriptstyle \mathcal{R}}\mathfrak{M}$
denote the category
of left $\mathcal R$-modules.

\begin{definition}
Two unital rings $\mathcal{R}$ and $\mathcal S$ are called
\emph{Morita equivalent} if
$_{\scriptscriptstyle \mathcal{R}}\mathfrak{M}$ and  
$_{\scriptscriptstyle \mathcal{S}}\mathfrak{M}$
are equivalent categories.
\end{definition}

\begin{example}\label{MatrixEx}
Let $\mathcal{R}$ be a unital ring and $M_n(\mathcal{R})$ be the
ring of $n \times n$ matrices over $\mathcal{R}$. Given a left
$\mathcal{R}$-module $V$,
we can define a left $M_n(\mathcal{R})$-module $\mathcal{F}(V)=V^n$, with the
$M_n(\mathcal{R})$-action given by matrix operating on vectors. The
functor $\mathcal{F}: {_{\scriptscriptstyle \mathcal{R}}\mathfrak{M}}\longrightarrow
{_{\scriptscriptstyle M_n(\mathcal{R})}\mathfrak{M}}$ defines an equivalence of
categories \cite[Thm.~17.20]{Lam99},
and $\mathcal{R}$ and $M_n(\mathcal{R})$ are Morita equivalent.
\end{example}

Isomorphic rings are Morita equivalent.
Properties of a ring $\mathcal{R}$ which are preserved
under Morita Equivalence are
called {\it Morita invariants}. Example \ref{MatrixEx} shows that commutativity
is {\it not} a Morita invariant property.
However, properties such as $\mathcal{R}$ being semisimple,
artinian and noetherian
are Morita invariant. Morita equivalent rings have isomorphic $K$-theory,
isomorphic lattice of ideals and isomorphic centers;
hence two commutative unital rings are Morita equivalent
if and only if they are isomorphic \cite[Cor.~18.42]{Lam99}.

Let  $\mathcal{R}$ and $\mathcal S$ be unital rings. An
$(\mathcal S,\mathcal R)$-bimodule $\SER$ canonically defines a functor
$\mathcal{F} = ({\SER} \otimes_{\scriptscriptstyle \mathcal{R}} \cdot) :
{_{\scriptscriptstyle \mathcal{R}}\mathfrak{M}} \longrightarrow
{_{\scriptscriptstyle \mathcal{S}}\mathfrak{M}}$  by
$$
\mathcal{F}({_{\scriptscriptstyle \mathcal{R}}V}) = 
{\SE}\otimes_{\scriptscriptstyle \mathcal{R}} V.
$$
It is  clear that $\mathcal{F}({_{\scriptscriptstyle \mathcal{R}}V})$
has a natural $\mathcal S$-module structure
determined by $s(x \otimes v) = sx \otimes v, \, x \in \mathcal{E}, \; v\in V$.
If $f:{_{\scriptscriptstyle \mathcal{R}}V}_1 \longrightarrow 
{_{\scriptscriptstyle \mathcal{R}}V}_2$
is a morphism, then we define
$\mathcal{F}(f): {\SE}\otimes_{\scriptscriptstyle \mathcal R}
V_1 \longrightarrow
{\SE}\otimes_{\scriptscriptstyle \mathcal R}V_2$
by $\mathcal{F}(f)(x\otimes v)= x\otimes f(v),\; x\in \mathcal{E}, \,
v\in V_1$.

It turns out that this way of constructing functors is very general.
It follows from a theorem of Eilenberg and Watts \cite{Watts60} that if
$\mathcal{F}: {_{\scriptscriptstyle \mathcal{R}}}\mathfrak{M} \longrightarrow
{_{\scriptscriptstyle \mathcal{S}}\mathfrak{M}}$
is an equivalence of categories, then there exists a bimodule $\SER$
such that $ {\SE} \otimes_{\scriptscriptstyle \mathcal R} \cdot \cong \mathcal{F}$. 
Under this  identification, the composition of functors corresponds to the 
tensor product of bimodules.

\begin{example}
In the case of $\mathcal{R}$ and $M_n(\mathcal{R})$, the functor described
in Example \ref{MatrixEx} corresponds to  the bimodule
$_{\scriptscriptstyle M_n(\mathcal{R})}\mathcal{R}^n_{\scriptscriptstyle \mathcal{R}}$.
\end{example}

\begin{corollary}
Two unital rings $\mathcal{R}$ and $\mathcal S$ are Morita equivalent
if and only
there exist bimodules $\SER$ and $\RES$
so that
$\RE\otimes_{\scriptscriptstyle \mathcal{S}} \ER \cong 
{_{\scriptscriptstyle \mathcal{R}}\mathcal{R}}_{\scriptscriptstyle \mathcal{R}}$ and
$\SE\otimes_{\scriptscriptstyle \mathcal{R}} \ES \cong 
{_{\scriptscriptstyle \mathcal{S}}\mathcal{S}}_{\scriptscriptstyle \mathcal{S}}$
(as bimodules).
\end{corollary}

\noindent A bimodule $\SER$ establishing a Morita equivalence is called
an \emph{equivalence bimodule}. Morita's theorem provides a characterization
of such bimodules.

\begin{definition}\label{ProgDef}
A right $\mathcal{R}$-module $\ER$ is
called a \emph{progenerator} if it is finitely generated,
projective and a generator\footnote{Recall that a right
$\mathcal{R}$-module $\ER$ is a {\it generator}
if any other right $\mathcal R$-module can be
obtained as a quotient of a direct sum of
copies of $\ER$.}.
\end{definition}

\begin{theorem}[Morita's theorem]
Two unital rings $\mathcal{R}$ and $\mathcal S$ are Morita equivalent if and only
if there exists a progenerator right $\mathcal{R}$-module $\ER$
so that $\mathcal{S} \cong \End_{\scriptscriptstyle \mathcal R}(\ER)$.
Moreover, if $\SER$ is an equivalence bimodule, then
its inverse is given by $\RES=\Hom_{\scriptscriptstyle \mathcal R}(\ER,\mathcal{R})$.
\end{theorem}
\begin{definition}\label{FullDef}
An idempotent $P\in M_n(\mathcal{R})$ ($P^2=P$) is called \emph{full} if
the span of elements of the form $T P S$, with $T,S \in M_n(\mathcal R)$,
is $M_n(\mathcal{R}) $. 
\end{definition}
A finitely generated projective $\mathcal{R}$-module
$\ER = P\mathcal{R}^n$ ($P\in M_n(\mathcal{R})$ idempotent)
is a generator if and only
if $P$ is full \cite[Remark 18.10(D)]{Lam99}. This provides an alternative description of
Morita equivalent rings.
\begin{theorem}\label{FullMorThm}
$\mathcal{R}$ and $\mathcal S$ are Morita equivalent
 if and only if there exists $n \in
\mathbb N$ and a full idempotent $P\in M_n(\mathcal{R})$ so that
$\mathcal{S}\cong PM_n(\mathcal{R})P = \End_{\scriptscriptstyle \mathcal{R}}
(P\mathcal{R}^n)$.
\end{theorem}

We note that there is a natural group associated to any unital ring $\mathcal R$.

\begin{definition}\label{PicDef}
We define $\Pic(\mathcal{R})$ as the group of
equivalence classes of self-equi\-va\-len\-ce
functors $\mathcal{F}:{_{\scriptscriptstyle \mathcal{R}}\mathfrak{M}}
\longrightarrow {_{\scriptscriptstyle \mathcal{R}}\mathfrak{M}}$, with group
operation given by composition; 
equivalently, we can view $\Pic(\mathcal{R})$ as the group of
isomorphism classes of
$(\mathcal{R},\mathcal{R})$-equivalence
bimodules $\RER$, with group operation
given by tensor products (over $\mathcal{R}$).
\end{definition}
\noindent The group $\Pic(\mathcal{R})$ is called the
\emph{Picard group} of $\mathcal{R}$.
\begin{remark}\label{PicRem}
Note that if $\RER$ is an equivalence bimodule, then the center of $\mathcal{R}$
need not act the same on the left and right of $\mathcal E$. If $\mathcal{R}$
is commutative and $\mathcal{E}$ is an $(\mathcal{R},\mathcal{R})$-equivalence
bimodule, then there exists an $(\mathcal{R},\mathcal{R})$-equivalence
bimodule $\mathcal{E}'$ satisfying $r x = x r$, for all $r \in \mathcal{R}$
and $x \in \mathcal{E}'$, such that $\mathcal{E} \cong \mathcal{E}'$ as
right $\mathcal{R}$-modules (pick $\mathcal{E}'$ of the form $P_0\mathcal{R}^n$
as a right $\mathcal{R}$-module
and consider the identification $\mathcal{R}\longrightarrow P_0M_n(\mathcal{R})P_0$,
$r \mapsto rP_0$).
\end{remark}

\noindent If $\mathcal{R}$ is commutative, we denote the group of isomorphism classes of 
$(\mathcal{R},\mathcal{R})$-equivalence bimodules $\RER$ satisfying 
$r x = x r$, for all $x \in \mathcal{E}$, $r \in \mathcal{R}$, by 
$\Pic_{\scriptscriptstyle{\mathcal{R}}}(\mathcal{R})$.

\begin{example}\label{PicMEx}
Let $\mathcal R = C^\infty(M)$, where $M$ is  a smooth manifold. As a consequence of
Serre-Swan's theorem
\cite[Chp.~XIV]{Bass68} (here used in the smooth category, where
the compactness assumption can be dropped),
$\Pic_{\scriptscriptstyle{C^\infty(M)}}(C^\infty(M))$ can be identified
with $\Pic(M)$, the group of isomorphism classes of complex line bundles
over $M$, with group operation given by fiberwise tensor product.
The Chern class map $\ch_1:\Pic(M) \longrightarrow H^2(M,\mathbb Z)$ is a
group isomorphism \cite[Sec.~3.8]{Hirz95}, and hence 
$\Pic_{\scriptscriptstyle{C^\infty(M)}}(C^\infty(M))\cong
H^2(M,\mathbb Z)$.
\end{example}

\subsection{Deformations of associative algebras: star products}\label{DefSub}
Let $k$ be a commutative, unital ring of characteristic zero,
and let $\mathcal A$ be an associative
$k$-algebra; $\mathcal{A}[[\lambda]]$ denotes the space of formal power series
with coefficients in $\mathcal A$.

\begin{definition}
A \emph{formal deformation} of $\mathcal A$ is an associative
$k[[\lambda]]$-bilinear multiplication $\star$ on
$\mathcal{A}[[\lambda]]$ of the form
\begin{equation}\label{FormalEq}
    A \star A' = \sum_{r=0}^\infty C_r(A,A') \lambda^r,
    \qquad A,A' \in \mathcal A,
\end{equation}
where the maps $C_r:\mathcal{A} \times \mathcal{A}\longrightarrow \mathcal A$ 
are $k$-bilinear, and $C_0$ is the
original product on $\mathcal A$ \footnote{We
extend $\star$ to $\mathcal{A}[[\lambda]]$ using $\lambda$-linearity.}.
\end{definition}
\noindent A formal deformation of $\mathcal A$ will be denoted by
$\qA=(\mathcal{A}[[\lambda]],\star)$.

\begin{definition}
We say that two 
formal deformations of $\mathcal A$, $\qA_1=(\mathcal{A}[[\lambda]],\star_1)$ and
$\qA_2=(\mathcal{A}[[\lambda]],\star_2)$, are  \emph{equivalent}
if there exist $k$-linear maps
$T_r: \mathcal{A}\longrightarrow \mathcal{A}$, $r \geq 1$, so that
$T =
\id + \sum_{r=1}^\infty T_r \lambda^r : \qA_1 \longrightarrow \qA_2$
satisfies
\begin{equation}
    A \star_1 A' = T^{-1}(T(A)\star_2 T(A')),
    \qquad \forall A, A' \in \mathcal A[[\lambda]].
\end{equation}
\noindent Such a $T$ is called an \emph{equivalence transformation}.
\end{definition}
\noindent The set of equivalence classes of formal deformations of $\mathcal A$
is denoted by $\Defa(\mathcal A)$. We denote the equivalence class of a deformation
$\star$ by $[\star]$.

If $\mathcal A$ is unital, then so is any formal deformation $\qA$; moreover,
any formal deformation of $\mathcal A$ is equivalent to one for which the unit
is the same as for $\mathcal A$ \cite[Sec.~14]{GS88}.

The group of automorphisms of $\mathcal{A}$,
 $\Aut(\mathcal{A})$, 
acts on formal deformations by $\star' = \psi^*(\star)$ if and only if
$A \star' A' = \psi^{-1}(\psi(A)\star \psi(A')), \; A, A' \in \mathcal{A}$.
Since any $k[[\lambda]]$-algebra isomorphism 
$S: (\mathcal{A}[[\lambda]],\star) 
\longrightarrow (\mathcal{A}[[\lambda]],\star')$ is of the form
$S= S_0 + \sum_{r=1}^\infty \lambda^r S_r$, with $k$-linear
$S_r:\mathcal{A} \longrightarrow \mathcal{A}$, and 
$S_0 \in \Aut(\mathcal{A})$,
a simple computation shows the following proposition.

\begin{proposition}\label{prop:iso}
Let $\star$ and $\star'$ be formal deformations of $\mathcal{A}$. Then
they are isomorphic if and only if there exists $\psi \in \Aut(\mathcal{A})$
with $[\psi^*(\star')] = [\star]$.
\end{proposition}

Let $\mathcal A$ be commutative and unital. 
\begin{definition}\label{BrackDef}
A \emph{Poisson bracket} on
$\mathcal A$ is a Lie algebra bracket $\{\cdot,\cdot\}$ satisfying the
Leibniz rule
$$
\{A_1,A_2A_3\} = \{A_1,A_2\} A_3 + A_2\{A_1,A_3\}.
$$
The pair $(\mathcal A,\{\cdot,\cdot\})$ is called a \emph{Poisson algebra}.
\end{definition}
\noindent For a formal deformation of $\mathcal A$ (\ref{FormalEq}), a simple
computation using associativity of $\star$ shows that
\begin{equation}\label{BrackEq}
\{A_1,A_2\}:= C_1(A_1,A_2)-C_1(A_2,A_1)= \frac{1}{\lambda}
(A_1\star A_2 - A_2\star A_1)\mbox{ mod } \lambda
\end{equation}
is a Poisson bracket on $\mathcal A$ \cite[Sect.~16]{SilWein99}.
It is simple to check that if two formal deformations are equivalent, then they
determine the same Poisson bracket through (\ref{BrackEq}).

Let $\mathcal A = C^\infty(M)$ be the algebra of complex-valued smooth functions
on a smooth manifold $M$. 

\begin{definition}\label{DefDef}
A formal deformation 
$\star = \sum_{r=0}^\infty C_r\lambda^r$ of $\mathcal A$ is called a 
\emph{star product} if each $C_r$ is a bidifferential operator. 
\end{definition}
\noindent The set of equivalence classes of star products on $M$ is denoted  by
$\Def(M)$.

A \emph{Poisson structure} on a manifold $M$ is a Poisson bracket $\{\cdot,\cdot\}$
on $C^\infty(M)$; it can be equivalently defined by a bivector field $\pi 
\in \Gamma^\infty(\bigwedge^2TM)$ satisfying $[\pi,\pi]=0$,\footnote{$[\cdot,\cdot]$
is the \emph{Schouten bracket} \cite{Vais94}.} in such a way that
$ \{f,g\} = \pi(df,dg)$. We call $\pi$ a \emph{Poisson tensor} and $(M,\pi)$
a \emph{Poisson manifold}.
We let
\begin{equation}\label{DefpiDef}
\Def(M,\pi):= \{[\star] \in \Def(M)\; | \; f\star g - g \star f = \lambda
\pi(df,dg) \mbox{ mod } \lambda^2\}.
\end{equation}

Star products on a Poisson manifold  $(M,\pi)$
will be assumed to be compatible with the
given Poisson structure in the sense of (\ref{BrackEq}).

If $\star$ and $\star'$ are star products on $(M,\pi)$,
we remark that, by Proposition \ref{prop:iso}, they
are isomorphic if and only if there exists
a Poisson diffeomorphism $\psi: M \longrightarrow M$ such that
$[\psi^*(\star')]=[\star]$. 

\subsection{Deformations of projective modules and vector bundles}\label{section:defvb}

Let $\mathcal A$ be a unital $k$-algebra, and let $\mathcal E$ be a
right  module over $\mathcal A$.
Let
$\Right_0: \mathcal{E} \times \mathcal{A} \longrightarrow \mathcal E$
denote the right action of $\mathcal A$ on $\mathcal E$,
$\Right_0(x, A) = x . A$.
Let $\qA = (\mathcal{A}[[\lambda]], \star)$ be a formal deformation of
$\mathcal A$ and suppose there exist $k$-bilinear maps
$\Right_r:\mathcal E \times \mathcal A \longrightarrow \mathcal E$,
for $r \geq 1$, such that the map
\begin{equation}
    \label{RrEq}
    \qR = \sum_{r=0}^\infty \Right_r \lambda^r :
    \mathcal{E}[[\lambda]] \times \qA
    \longrightarrow \mathcal{E}[[\lambda]]
\end{equation}
makes $\mathcal{E}[[\lambda]]$ into a module over $\qA$. We will
write $\qR(x,A) = x \bullet A$, for $x \in \mathcal E$,
$A \in \mathcal A$.
\begin{definition}\label{qEDef}
We call $\qE=(\mathcal{E}[[\lambda]], \bullet)$ a \emph{deformation}
of the right $\mathcal A$-module $\mathcal E$ corresponding to
$\qA=(\mathcal{A}[[\lambda]], \star)$. Two deformations
$\qE=(\mathcal{E}[[\lambda]], \bullet)$,
$\qE'=(\mathcal{E}[[\lambda]], \bullet')$ are \emph{equivalent} if
there exists an $\qA$-module isomorphism
$T: \qE \longrightarrow \qE'$ of the form
$T = \id + \sum_{r=1}^\infty T_r \lambda^r$,
with $k$-linear maps $T_r:\mathcal E \longrightarrow \mathcal E$.
\end{definition}

Suppose $\mathcal E$ is finitely generated and projective 
over $\mathcal A$, in which case we can write  $\mathcal E=P_0\mathcal{A}^n$
for some idempotent $P_0 \in M_n(\mathcal A)$ and $n \in \mathbb N$. 
Let $\qA$ be a formal deformation
of $\mathcal A$. It is clear that  $M_n(\qA)$ can be identified with
$M_n(\mathcal A)[[\lambda]]$ as a $k[[\lambda]]$-module, and this
identification naturally
defines a formal deformation of $M_n(\mathcal A)$.
We recall the following well-known fact \cite{Fed96,EmmWe96,GS88}.

\begin{lemma}\label{ProjDefLem}
Let $P_0 \in M_n(\mathcal A)$ be
an idempotent. For any formal deformation
$\qA$, there exists an idempotent $\qP \in M_n(\qA)$ with
$\qP = P_0 + O(\lambda)$. 
\end{lemma}

\noindent We call $\qP$ a \emph{deformation} of $P_0$ with respect to $\qA$.
The natural right $\qA$-module structure on $\qP \star \qA^n$
can then be transfered to $P_0\mathcal{A}^n[[\lambda]]$ 
since $P_0\mathcal{A}^n[[\lambda]]
\cong \qP \star \qA^n$ as $k[[\lambda]]$-modules. An explicit isomorphism
is given by the map \cite[Lem.~2.3]{BuWa2000b}
\begin{equation}\label{IEq1}
J:P_0\mathcal{A}^n[[\lambda]]\longrightarrow \qP \star \qA^n, \;\; P_0x \mapsto
\qP\star x,
\end{equation}
$x \in \mathcal{A}^n$. By Lemma \ref{ProjDefLem}, deformations of finitely generated
projective modules (f.g.p.m.) always exist
(with respect to any $\qA$). Moreover, they are
unique up to equivalence \cite[Prop.~2.6]{BuWa2000} (and hence necessarily
finitely generated and projective over $\qA$).
We summarize these facts in the following proposition.

\begin{proposition}\label{prop:unique}
Let $\EA$ be a finitely generated projective module over $\mathcal A$,
and let $\qA=(\mathcal{A}[[\lambda]],\star)$ be a formal deformation.
Then there exists a deformation of $\EA$ with respect to $\qA$, and any
two such deformations are equivalent.
\end{proposition}

A simple computation shows that fullness of 
idempotents is preserved under deformations:

\begin{lemma}\label{DefFullLem}
Let $P_0 \in M_n(\mathcal A)$ be an idempotent and $\qP \in M_n(\qA)$
be a deformation of $P_0$. Then $P_0$ is full if and only if $\qP$ is full.
\end{lemma}

\begin{proposition}\label{BijecProp}
An equivalence bimodule
 $\BEA$ canonically defines
a bijective map  $\Phi_{\scriptscriptstyle{\mathcal{E}}}:
\Defa(\mathcal{A}) \longrightarrow \Defa(\mathcal{B})$ so that
formal deformations related by $\Phi$ are Morita equivalent.
\end{proposition}

\begin{proof}
Let $\star$ be a formal deformation of $\mathcal{A}$, and let
$\qE=(\mathcal{E}[[\lambda]],\bullet)$ be a deformation of $\EA$
with respect to $\star$. Then $\mathcal{B}[[\lambda]]$ and 
$\End(\qE)$ are isomorphic as $k[[\lambda]]$-modules. In fact,
for $\mathcal{E}=P_0\mathcal{A}^n$, $P_0$ full idempotent, and
$\qE$ given by (\ref{IEq1}), an explicit isomorphism is
(see \cite[Lem.~2.3]{BuWa2000b})
\begin{equation}\label{IEq2}
I: P_0M_n(\mathcal A)P_0[[\lambda]] \longrightarrow \qP \star M_n(\qA)  \star \qP,
\;\; P_0LP_0 \mapsto \qP\star L \star \qP.
\end{equation}
In general,
any fixed isomorphism 
\begin{equation}\label{IEq3}
\mathcal{B}[[\lambda]] \stackrel{\sim}{\longrightarrow}\End(\qE)
\end{equation} 
induces a formal
deformation $\star'$ 
of $\mathcal{B}$, and different choices of (\ref{IEq3})
 lead to isomorphic deformations (not necessarily equivalent). 
Note that an isomorphism (\ref{IEq3}) also defines
a left module structure $\bullet'$ on $\mathcal{E}[[\lambda]]$
over $(\mathcal{B}[[\lambda]],\star')$, and we can choose
(\ref{IEq3}) so that
$\bullet'$ is a deformation of
the module structure of $\BE$, $B \bullet' x = B\cdot x + O(\lambda)$ 
(this is the case for (\ref{IEq2})). A simple computation shows that
any two isomorphisms (\ref{IEq3}) yielding deformations of
$\BE$ define equivalent deformations of $\mathcal{B}$. Moreover,
it follows from the uniqueness part of
 Proposition \ref{prop:unique} that this equivalence class
does not depend on the choice of $\bullet$ or element in $[\star]
\in \Defa(\mathcal{A})$. Hence this procedure defines a canonical map
\begin{equation}
\Phi_{\scriptscriptstyle{\mathcal{E}}}: 
\Defa(\mathcal A) \longrightarrow \Defa(\mathcal B),
\end{equation}
which  is  a bijection \cite[Prop.~3.3]{BuWa2000b}. 

It follows from the construction of 
$\Phi_{\scriptscriptstyle{\mathcal{E}}}$ and Lemma \ref{DefFullLem}
that formal deformations related by 
$\Phi_{\scriptscriptstyle{\mathcal{E}}}$ are Morita equivalent.
\end{proof}

Let $\mathcal A = C^\infty(M)$ be the algebra
of complex-valued smooth
functions on a manifold $M$; by Serre-Swan's theorem, f.g.p.m. over $\mathcal A$
correspond to  (smooth sections of)
finite dimensional complex vector bundles $E \to M$.

For a smooth complex $m$-dimensional vector bundle
$E \to M$,  let $\mathcal E = \Gamma^\infty(E)$ be the space of smooth sections
of $E$,
regarded as a  right $\mathcal A$-module, and
$\mathcal B = \End_{\mathcal A}(\mathcal E) = \Gamma^\infty(\End(E))$, the
algebra of smooth sections of the endomorphism bundle $\End(E) \to M$.

Let $\star$ be a star product on $M$.

\begin{definition}
\label{DeqQuanDef}
A \emph{deformation quantization} of $E$ with respect to $\star$
is a deformation of $\mathcal E = \Gamma^\infty(E)$ in the sense of
Definition \ref{qEDef} so that the
corresponding $\Right_r$ (as in (\ref{RrEq}))
are bidifferential operators.
\end{definition}

\noindent The explicit map in (\ref{IEq1}) shows that the (always existing)
deformations of $\mathcal E$ can be chosen with bidifferential $\Right_r$
(since we are only considering differential star products); therefore
deformation quantizations of vector bundles exist (with respect to any star product)
 and are unique up to equivalence.

If we write $\mathcal E = P_0 \mathcal{A}^n$, for some idempotent 
$P_0 \in M_n(\mathcal A)$, we note that the deformation of $\mathcal B$
induced by the explicit map (\ref{IEq2}) as in (\ref{IEq3}) is also differential.
This gives rise to a canonical bijective map (Proposition \ref{BijecProp})
\begin{equation}\label{BijecDiffEq}
\Phi_{\scriptscriptstyle \mathcal E}: 
\Def(M) \longrightarrow  \Def(\Gamma^\infty(\End(E))).
\end{equation}

\subsection{Poisson fibred algebras}\label{PoissonFibredSub}
 Poisson fibred algebras \cite{RVW96} arise in connection with formal deformations of
 noncommutative algebras, providing a generalization of Poisson algebras. We will 
recall the definition.

Let $\mathcal B$ be a unital $k$-algebra, not necessarily commutative, and
let $\mathcal Z$ be its center. 
\begin{definition}
A \emph{Poisson fibred algebra structure} on $(\mathcal B, \mathcal Z)$ is
a bracket
$$
\{\cdot,\cdot\} : \mathcal{Z} \times \mathcal{B} \longrightarrow \mathcal{B}
$$
satisfying the following conditions:
\begin{enumerate}
\item[1.] The restriction of $\{\cdot,\cdot\}$ to $\mathcal{Z} \times \mathcal{Z}$
makes $\mathcal Z$ into a Poisson algebra.

\item[2.] The following Leibniz identities hold \footnote{These identities imply that
$\{Z,1\}=\{1,B\}=0$, for all $Z \in \mathcal Z,\; B \in \mathcal B$.}.
\begin{gather}
\{Z, B_1 B_2\} = \{Z,B_1\}B_2 + B_1\{Z,B_2\},\label{FibredEq1}\\
\{Z_1 Z_2, B\} = Z_1\{Z_2,B\} + Z_2\{Z_1,B\}. \label{FibredEq2}
\end{gather}
\end{enumerate}
\end{definition}

Suppose  $\mathcal A$ is a commutative, unital $k$-algebra, 
and let $P_0 \in M_n(\mathcal A)$
be a full idempotent. Let $\qA=(\mathcal{A}[[\lambda]],\star)$ be a formal
deformation of $\mathcal A$. We keep the notation
$$
A_1 \star A_2 = \sum_{r=0}^\infty C_r(A_1,A_2)\lambda^r,\;\; A_1,A_2 \in
\mathcal A.
$$
\noindent Since $\mathcal A$ is commutative, it inherits a Poisson
algebra structure  (\ref{BrackEq}) from $\star$ given by
$$
\{A_1,A_2\}:=C_1(A_1,A_2)-C_1(A_2,A_1).
$$

We saw in Section \ref{section:defvb} how to define a formal deformation $\star'$
of $\mathcal B = P_0 M_n(\mathcal A)P_0$ explicitly by choosing an idempotent
$\qP \in M_n(\qA)$ deforming $P_0$: 
$$
L_0\star'S_0 := I^{-1}(I(L_0)\star I(S_0))=
\sum_{r=0}^\infty B_r(L_0,S_0)\lambda^r, \;\; L_0,S_0 \in \mathcal B,
$$
where $I$ is as in (\ref{IEq2}).
For $M,N \in M_n(\mathcal A)$, let (note the abuse of notation\footnote{The
bracket in $M_n(\mathcal A)$ defined in (\ref{MatrixBrackEq})
induces the bracket $\{\cdot,\cdot\}$ on $\mathcal{A}$ through the identification
of $\mathcal A$ with the center of $M_n(\mathcal A)$ in the natural way;
we denote both brackets by $\{\cdot,\cdot\}$ by abuse of notation.})
\begin{equation}\label{MatrixBrackEq}
\{M,N\}:= C_1(M,N)-C_1(N,M),
\end{equation}
where $C_1(M,N)\in M_n(\mathcal A)$ is defined by
${C_1(M,N)}_{i,j}= \sum_{r=1}^nC_1(M_{i,r},N_{r,j})$.
We will compute the expression for the bracket in $P_0M_n(\mathcal A)P_0
\subseteq M_n(\mathcal A)$ given by
\begin{equation}
\{L,S\}':= B_1(L,S)-B_1(S,L), \;\; L,S \in \mathcal{B},
\end{equation}
in terms of $\{\cdot,\cdot\}$.
\begin{proposition}\label{FibredProp}
For all $L_0,S_0 \in P_0M_n(\mathcal A)P_0$, we have
$\{L_0,S_0\}'=P_0\{L_0,S_0\}P_0.$
\end{proposition}
\begin{proof}
Since $I:P_0M_n(\mathcal A)P_0[[\lambda]] \longrightarrow
\qP \star M_n(\qA) \star \qP$ and 
$\qP \star M_n(\qA) \star \qP \subseteq M_n(\qA)=M_n(\mathcal A)[[\lambda]]$,
we can write $I=\sum_{r=0}^\infty I_r\lambda^r$, where
$$
I_r:P_0M_n(\mathcal A)P_0 \longrightarrow M_n(\mathcal A).
$$
A simple computation shows that $I_0(L_0)=L_0$ and
\begin{equation}\label{I1Eq}
I_1(L_0)=C_1(P_0,L_0) + P_0C_1(L_0,P_0) + L_0P_1 + P_1L_0.
\end{equation}
\noindent The equations $\qP\star I(L_0)=I(L_0)$ and $I(L_0)\star\qP = I(L_0)$
imply that
\begin{equation}\label{2Eq}
P_0C_1(P_0,L_0) + P_0P_1L_0 =0, \; \mbox{ and } \;
C_1(L_0,P_0)P_0 + L_0P_1P_0=0.
\end{equation}
\noindent It follows from (\ref{I1Eq}) and (\ref{2Eq}) that
$P_0I_1(L_0)P_0 =0$, for all $L_0\in P_0M_n(\mathcal A)P_0$.
Note that if $L_0+\lambda L_1 + \ldots \in P_0M_n(\mathcal A)P_0[[\lambda]]$
and $I(L_0+\lambda L_1 + \ldots) = M_0 + \lambda M_1 + \ldots$, then
$I_1(L_0) + I_0(L_1) = I_1(L_0) + L_1 =M_1$, and hence
$$
P_0(I_1(L_0)+L_1)P_0 = L_1 = P_0M_1P_0.
$$
But $I(L_0 \star' S_0)= I(L_0)\star I(S_0)= L_0S_0 +
\lambda (C_1(L_0,S_0)+L_0I_1(S_0)+I_1(L_0)S_0) + \dots$.
Thus
$$
B_1(L_0,S_0)=P_0(C_1(L_0,S_0)+L_0I_1(S_0)+I_1(L_0)S_0)P_0=P_0C_1(L_0,S_0)P_0,
$$
by (\ref{2Eq}), and the result follows.
\end{proof}

 Let $\mathcal Z$ denote the center of $\mathcal B = P_0M_n(\mathcal A)P_0$,
 $P_0 \in M_n(\mathcal A)$ full idempotent.
The triple  $(\mathcal B, \mathcal Z, \{\cdot,\cdot\}')$ is a Poisson fibred
algebra \cite[Prop.~1.2]{RVW96}, and, as such, the restriction
of $\{\cdot,\cdot\}'$ to $\mathcal Z \times \mathcal Z$ provides $\mathcal Z$
with the structure of a Poisson algebra.
We can identify $\mathcal Z$ and $\mathcal A$ through the 
algebra isomorphism
$\Psi:\mathcal A \longrightarrow \mathcal Z$, $A \mapsto P_0AP_0=AP_0$.

\begin{theorem}\label{BracketThm}
The map $\Psi: (\mathcal A, \{\cdot,\cdot\})\longrightarrow
(\mathcal Z,\{\cdot,\cdot\}')$ is an isomorphism of Poisson algebras.
\end{theorem}
\begin{proof}
The bracket $\{\cdot,\cdot\}':\mathcal{Z} \times \mathcal B \longrightarrow
\mathcal B$ satisfies (\ref{FibredEq1}) and (\ref{FibredEq2}).
\noindent As a result,
$\{\cdot,P_0\}'=P_0\{\cdot,P_0\}P_0 = 0$ and $\{P_0,\cdot\}'=
P_0\{P_0,\cdot\}P_0=0$. It is easy to check that the following
Leibniz rule holds  for
$\{\cdot,\cdot\}$  in $M_n(\mathcal A)$:
$$
\{A M, N\}=A\{M,N\} + M\{A,N\}, \; M,N \in M_n(\mathcal A),\;\,
A \in \mathcal A \cong \mbox{center}(M_n(\mathcal A)).
$$
Combining these identities, we get
$$
\{\Psi(A_1),\Psi(A_2)\}'=\{A_1P_0,A_2P_0\}'=\{A_1,A_2\}P_0=\Psi(\{A_1,A_2\}),
$$
for $A_1,A_2 \in \mathcal A$.
\end{proof}

\section{Picard groups acting on deformations}\label{ActionSec}

Let $\mathcal A$ be a commutative, unital $k$-algebra. 

By Proposition \ref{BijecProp}, 
an $(\mathcal{A},\mathcal{A})$-equivalence
bimodule $\AEA$ canonically defines
an automorphism of the set $\Defa(\mathcal A)$,
\begin{equation}\label{AutoEq}
\Phi_{\scriptscriptstyle \mathcal E}:
\Defa(\mathcal A) \longrightarrow \Defa(\mathcal A).
\end{equation}
We observe that the map $\Phi_{\scriptscriptstyle \mathcal E}$ only
depends on the isomorphism class of $\mathcal{E}$. We will abuse notation
and simply write $\mathcal{E}$ to denote its isomorphism class
in $\Pic(\mathcal A)$ or $\Pic_{\scriptscriptstyle{\mathcal{A}}}(\mathcal{A})$
(if $Ax = x A$ for all $x \in \mathcal{E}$, $A \in \mathcal{A}$).

Let $\qA=(\mathcal{A}[[\lambda]],\star),
\qA'=(\mathcal{A}[[\lambda]],\star')$ be formal deformations of $\mathcal A$.

\begin{lemma}\label{MoritaActlem} The unital $k[[\lambda]]$-algebras
$\qA$ and $\qA'$ are Morita equivalent if and only if there exists
$\mathcal E \in \Pic_{\scriptscriptstyle{\mathcal{A}}}(\mathcal A)$ and 
$\psi \in \Aut(\mathcal A)$ with 
$\Phi_{\scriptscriptstyle \mathcal E}([\star])=[\psi^*(\star')]$.
\end{lemma}
\begin{proof}
If $\Phi_{\scriptscriptstyle \mathcal E}([\star])=[\psi^*(\star')]$, 
then $\qA$ and $\qA'$ are Morita
equivalent by Propositions \ref{prop:iso} and  \ref{BijecProp}.
Conversely, if $\qA$ and $\qA'$ are Morita equivalent, then
there exists a full idempotent $\qP \in M_n(\qA)$ so that
$\qA' \cong \qP\star M_n(\qA)\star \qP$.
By Lemma \ref{DefFullLem}, $\qP = P_0 + O(\lambda)$ with $P_0$ full.
We know that (see (\ref{IEq1})) $\qP\star M_n(\qA)\star \qP$ is
isomorphic to
$P_0M_n(\mathcal{A})P_0[[\lambda]]$ as a $k[[\lambda]]$-module, and since it is
also isomorphic to $\mathcal A[[\lambda]]$, we must have
$P_0M_n(\mathcal{A})P_0 \cong \mathcal A$. As in Remark \ref{PicRem},
  $\mathcal E =
P_0\mathcal{A}^n$ is 
an $(\mathcal{A},\mathcal{A})$-equivalence bimodule
satisfying $xA=Ax$, for all $x \in \mathcal{E}$ and $A \in \mathcal{A}$, 
and $\qA'$ is isomorphic
to the deformations in the class 
$\Phi_{\scriptscriptstyle \mathcal E}([\star])$. The result then follows
from Proposition \ref{prop:iso}.
\end{proof}

\noindent We recall that the unit element in 
$\Pic_{\scriptscriptstyle{\mathcal{A}}}(\mathcal A)$ is given
by (the isomorphism class of) $\AAA$.

\begin{lemma} \label{Action1Lem}
If $\AEA \cong \AAA$, then $\Phi_{\scriptscriptstyle \mathcal{E}} = \id$.
\end{lemma}
\begin{proof}
Let $\qA=(\mathcal{A}[[\lambda]],\star)$ be a formal deformation of
$\mathcal A$. Then $\qA$ itself, regarded as a right module over $\qA$,
provides a deformation of $\EA$. Since $\End_{\qA}(\qA) = \qA$, it follows that
$\Phi_{\scriptscriptstyle \mathcal E}([\star])=[\star]$.
\end{proof}
\begin{lemma}
Let $\mathcal E$, $\mathcal E' \in 
\Pic_{\scriptscriptstyle{\mathcal{A}}}(\mathcal A)$, and $\mathcal{E}''=
\mathcal{E}'\otimes_{\mathcal A} \mathcal E$. Then
$\Phi_{\scriptscriptstyle{\mathcal{E}''}}=\Phi_{\scriptscriptstyle{\mathcal{E}'}}
\circ \Phi_{\scriptscriptstyle \mathcal E}$.
\end{lemma}
\begin{proof}
Let $\qA=(\mathcal{A}[[\lambda]],\star)$,
$\qA'=(\mathcal{A}[[\lambda]],\star')$, and
$\qA''=(\mathcal{A}[[\lambda]],\star'')$ be formal deformations of
$\mathcal A$ so that
$[\star']=\Phi_{\scriptscriptstyle \mathcal E}([\star])$ and
$[\star'']=\Phi_{\scriptscriptstyle{\mathcal E'}}([\star'])$.
Let $\qE$ be a deformation of $\mathcal E$ corresponding to $\star$,
and $\qE'$ be a deformation corresponding to $\star'$ .
We know that
$\qA' \cong \End_{\qA}(\qE)$ and $\qA''\cong \End_{\qA'}(\qE')$.
As discussed in Section \ref{MorPicSub}, $\qE'\otimes_{\qA'} \qE$ is
an $(\qA'',\qA')$-equivalence bimodule, so
$\qA'' \cong \End_{\qA}(\qE'\otimes_{\qA'} \qE)$. Since
$\qE'\otimes_{\qA'} \qE$ is a f.g.p.m. over $\qA$, it follows
(see (\ref{IEq1})) that it is of the form
$V[[\lambda]]$, where $V \cong \qE'\otimes_{\qA'} \qE/
(\lambda \qE'\otimes_{\qA'} \qE)$ as a $k$-module.
But 
$$
\qE'\otimes_{\qA'} \qE/(\lambda \qE'\otimes_{\qA'} \qE) \cong
\mathcal{E}'\otimes_{\mathcal A}\mathcal{E}.
$$
Hence $\qE'\otimes_{\qA'} \qE$ is a deformation  of
$\mathcal{E}'\otimes_{\mathcal A} \mathcal{E}$, and the conclusion follows.
\end{proof}

\noindent The following lemma is a simple corollary of Theorem \ref{BracketThm}.

\begin{lemma}\label{PresProp}
Let $\star = \sum_{r=0}^\infty C_r \lambda^r$ and
$\star' = \sum_{r=0}^\infty C_r' \lambda^r$ be formal deformations
of $\mathcal A$ such that $[\star'] = \Phi_{\scriptscriptstyle \mathcal E}([\star])$,
for some $\mathcal E \in 
\Pic_{\scriptscriptstyle{\mathcal{A}}}(\mathcal A)$. Then $\star$ and $\star'$
correspond to the same Poisson bracket, i.e.
$$
C_1(A_1,A_2)-C_1(A_2,A_1) = C_1'(A_1,A_2)-C_1'(A_2,A_1),\; \mbox{ for all }
A_1,A_2 \in \mathcal{A}.
$$
\end{lemma}
\noindent As a result, we can state the following theorem.
\begin{theorem}\label{ActionThm}
Let $\mathcal A$ be a commutative, unital $k$-algebra.
Then $\Phi: \Pic_{\scriptscriptstyle{\mathcal{A}}}(\mathcal A) 
\to \Aut(\Defa(\mathcal{A}))$,
$\mathcal E \mapsto \Phi_{\scriptscriptstyle \mathcal E}$, defines an action
of $\Pic_{\scriptscriptstyle{\mathcal A}}(\mathcal A)$ 
on the set $\Defa(\mathcal A)$, preserving
Poisson brackets. Moreover, two
formal deformations of $\mathcal A$, $\star$ and $\star'$,
are Morita equivalent if and only if there exists 
$\psi \in \Aut(\mathcal{A})$ such that $[\star]$ and $[\psi^*(\star')]
$ lie in the same orbit of $\Phi$.
\end{theorem}

\section{Semiclassical geometry of quantum line bundles}\label{SemiClassSec}

Henceforth $\mathcal A = C^\infty(M)$, 
the algebra of complex-valued smooth functions on a manifold $M$,
and we will restrict ourselves
to differential deformations of $\mathcal A$ (i.e., star products).
As we noted in Example \ref{PicMEx}, $\Pic(M)\cong H^2(M,\mathbb Z)$
(group of isomorphism classes of complex line bundles over $M$) 
can be naturally identified with
$\Pic_{\scriptscriptstyle{\mathcal{A}}}(\mathcal A)$ 
through the map $\Pic(M) \ni L \mapsto \mathcal E = 
\Gamma^\infty(L) \in
\Pic_{\scriptscriptstyle{\mathcal{A}}}(\mathcal A)$.

The following result is a consequence of (\ref{BijecDiffEq}) 
and Theorem \ref{ActionThm}.
\begin{theorem}\label{StarActThm}
Let $(M,\pi)$ be a Poisson manifold.
There is a canonical action
$$
\Phi:\Pic(M)\times \Def(M,\pi) \longrightarrow \Def(M,\pi),
$$
and two star products $\star$ and $\star'$ on $(M,\pi)$ 
are Morita equivalent if and only if there exists a Poisson diffeomorphism
$\psi:M\longrightarrow M$ such that $[\star]$ and $[\psi^*(\star')]$ lie
in the same orbit of $\Phi$.
\end{theorem}

 The goal of the remainder of this work
is to understand the action $\Phi$ and the orbit space
$\Def(M,\pi)/\Pic(M).$ To this end, we will study the semiclassical
geometry of deformed line bundles over Poisson manifolds. A comparison
between the objects arising as ``first-order'' approximations to
deformed vector bundles and the usual notion of Poisson module \cite{RVW96}
is discussed in \cite{Bu2001b}.

Let $(M,\pi)$ be a Poisson manifold, and let
$\star = \sum_{r=0}^\infty C_r\lambda^r$
be a star product on $M$ satisfying
$$
C_1(f,g)-C_1(g,f)=\pi(df,dg),\; \,f,g \in C^\infty(M).
$$
Let $L \to M$ be a complex line bundle over $M$, and let
$\mathcal E = \Gamma^\infty(L)$.
Let us fix a deformation quantization of $L$ with respect to $\star$,
$\qE=(\mathcal{E}[[\lambda]],\bullet)$, and pick
$\star' \in \Phi_{\scriptscriptstyle \mathcal E}([\star])$.
Since $\qA'=(C^\infty(M)[[\lambda]],\star')
\cong \End_{\qA}(\qE)$, there is a natural left action of $\qA'$ on $\qE$ that can
be written
$$
f\bullet' s = fs + \sum_{r=1}^\infty R_r'(f,s),
$$
for bidifferential maps
$R_r':C^\infty(M) \times \mathcal{E} \longrightarrow \mathcal{E}$.
It is clear that $\bullet'$, $\bullet$ make $\mathcal{E}[[\lambda]]$ into an
$(\qA',\qA)$-bimodule.

\begin{definition}\label{QuantBimDef}
Let $L\to M$ be a complex line bundle over a Poisson manifold $M$.
Fix $\star=\sum_rC_r\lambda^r$ on $M$, and $\star'=\sum_rC_r'\lambda^r
\in \Phi_{\scriptscriptstyle \mathcal E}([\star])$.
A triple $(\mathcal{L}[[\lambda]],\bullet,\bullet')$
is called a bimodule quantization of $L$ corresponding to $\star$, $\star'$.
\end{definition}

The following equations relate $\star$, $\star'$, $\bullet$
and $\bullet'$.
\begin{eqnarray}
(f\star' g) \bullet s & = & f \bullet'(g \bullet' s), \\ \label{RelEqI}
s \bullet(f \star g) & = & (s \bullet f) \bullet g, \\ \label{RelEqII}
(f\bullet' s)\bullet g & = & f\bullet'(s\bullet g). \label{RelEqIII}
\end{eqnarray}

Let $R: \mathcal{E}\times C^\infty(M) \longrightarrow \mathcal{E}$ be defined
by
\begin{equation}
R(s,f):=R_1(s,f)-R_1'(f,s).
\end{equation}
Since $[\star]$ and $\Phi_{\scriptscriptstyle \mathcal E}([\star])$ 
correspond to the Poisson
bracket $\{\cdot,\cdot\}$ on $M$ (by Lemma \ref{PresProp}), we may
assume $C_1=C_1'$ in Definition \ref{QuantBimDef}.

\begin{proposition} \label{ConnProp} The map
$R$ is a contravariant connection on $L$.
\end{proposition}
\begin{proof}
We must show that $R$ satisfies (\ref{ContConnEqI}), (\ref{ContConnEqII}) in Appendix
\ref{AppA}.
Note that
(\ref{RelEqI})
yields, in first order,
\begin{equation}\label{I}
R_1'(fg,s)+ C_1'(f,g)s = R_1'(f,gs) + f R_1'(g,s).
\end{equation}
Similarly, (\ref{RelEqII})
yields
\begin{equation}\label{II}
R_1(s,fg)+s C_1(f,g) = R_1(sf,g)+R_1(s,g)f.
\end{equation}
We finally note that
(\ref{RelEqIII}) implies that
\begin{equation}\label{III}
R_1(fs,g)+R_1'(f,s)g=R_1'(f,sg)+fR_1(s,g).
\end{equation}
The difference of equations (\ref{I}) and (\ref{II}) yields
$$
 R(fg,s)=R_1(sf,g) + R_1(s,f)g -R_1'(f,gs) - fR_1'(g,s).
$$
But, by (\ref{III}), $R_1(sf,g)=R_1'(f,sg)+fR_1(s,g)-R_1'(f,s)g$.
This implies that
$$
R(s,fg)=fR(s,g)+ gR(s,f),
$$
proving that (\ref{ContConnEqI}) is satisfied.
Now, switching $f$ and $g$ in (\ref{I}), and subtracting it from (\ref{II})
(assuming $C_1=C_1'$),
we get
$$
R(sf, g) =  \{f,g\}s + R(s,fg) - gR(s,f) = \{f,g\}s +fR(s,g),
$$
proving (\ref{ContConnEqII}).
\end{proof}

We observe that given $\star$ on $M$, the contravariant connection
$R$ on $L$ depends on the choice of $\star'$, $\bullet$ and $\bullet'$.
As an example, let us compute it in a concrete situation.
\begin{example}\label{ConnEx}
Fix $n$, and
let $t(\mathbb{C}^n) = M \times \mathbb{C}^n \to M$ be a trivial bundle.
Let $P_0 \in M_n(C^\infty(M))$ be a rank-one idempotent
so that $L=P_0 t(\mathbb{C}^n)$ is a line bundle over $M$.
For a fixed star product $\star$ on $M$,
we pick an idempotent $\qP = P_0 + O(\lambda) \in M_n(\qA)$.
Using $J$ in (\ref{IEq1}), and $I$ in (\ref{IEq2}) to 
establish  $\mathbb{C}[[\lambda]]$-module
isomorphisms $P_0\mathcal{A}^n[[\lambda]] \longrightarrow
\qP\star \qA^n$, and $P_0M_n(\mathcal A)P_0[[\lambda]] \longrightarrow
\qP\star M_n(\qA)\star \qP$, respectively, an explicit computation (in the spirit
of Proposition \ref{FibredProp}) shows that
$R_1(s,f)= P_0C_1(s,f)$ and $R_1'(f,s)=P_0C_1(f,s)$, where
$C_1(f,s)_i=C_1(s_i,f)$ and $C_1(s,f)_i=C_1(f,s_i)$, $i=1,\ldots, n$.
Thus
$$
R(s,f) = P_0\{s,f\} = \nabla_{X_f}s,
$$
where $\nabla = P_0 d$ is the adapted
connection on $L$, and $X_f$ the hamiltonian vector field of $f$.
\end{example}

Let $D$ be a contravariant connection on $L$ induced by an ordinary connection
$\nabla$ (i.e., $D_{df}s=\nabla_{X_f}s$).
Fix $\star$ on $M$. 
It follows from Example \ref{ConnEx} and \cite[Thm.~1.1]{PoRe86} that
we can choose a bimodule quantization of $L$ so that $R=D$.
This in fact holds for any contravariant connection (not necessarily
induced by an ordinary one), as
 discussed in \cite{Bu2001b}.

\section{The semiclassical limit of Morita equivalent star products} \label{EquivSec}
\subsection{Semiclassical curvature}
Let $(M,\pi)$ be a Poisson manifold, and suppose
$\star=\sum_{r=0}^\infty C_r\lambda^r$  and
$\star'=\sum_{r=0}^\infty C_r'\lambda^r$ are star products on $M$, with
$C_1=C_1'=\frac{1}{2}\{\cdot,\cdot\}$. We can associate a Poisson cohomology
class to the pair $[\star],[\star']$, measuring the obstruction for
these star products being equivalent  modulo $\lambda^3$ \cite[Prop.~3.1]{BCG97}.
\begin{lemma}\label{TauLem}
Suppose $\star$ and $\star'$ are star products with
$C_1=C_1'=\frac{1}{2}\{\cdot,\cdot\}$.
The map
$$
(df,dg)\mapsto (C_2-C_2')(f,g) -(C_2-C_2')(g,f)
$$
defines a $d_\pi$-closed bivector field $\tau \in \Gamma^\infty(\bigwedge^2TM)$.
Moreover, the class $[\tau]_\pi \in H_\pi^2(M)$
depends only on the equivalence classes of $\star$ and $\star'$.
\end{lemma}
\begin{proof}
The fact that $\tau$ is a closed bivector field was proven in
\cite[Prop.~3.1]{BCG97}.

Suppose $\hat{\star}$ is a star product equivalent to $\star$:
$$
f \hat{\star} \, g = \sum_{r=0}^\infty \widehat{C}_r \lambda^r =
\qT^{-1}(\qT(f)\star \qT(g)),
$$
where $\qT = \id + \sum_{r=1}^\infty T_r\lambda^r$ is an equivalence
transformation. Assume $\widehat{C}_1=C_1=\frac{1}{2}\{\cdot,\cdot\}$.
We must
show that, if $\widehat{\tau}$ is the closed bivector given by the skew
symmetric part of $(\widehat{C}_2-C_2')$, then $[\widehat{\tau}]_\pi = [\tau]_\pi$.
The condition $\widehat{C}_1=C_1$ implies that $T_1 \in \Der(C^\infty(M))$.
Hence $T_1=\mathcal{L}_X$, for some vector field $X \in \chi(M)$.
A simple computation just using the definitions shows that
$$
\widehat{\tau}=\tau - d_\pi X,
$$
and the result follows.
\end{proof}
 
 If $M$ is a symplectic manifold, then the bivector
$\tau$ defines a closed $2$-form $\widetilde{\tau}$ by
\begin{equation}\label{TauTilEq}
\widetilde{\tau}(X_f,X_g)=\tau(df,dg),
\end{equation}
where $X_f$ and $X_g$ are the hamiltonian vector fields of $f$ and $g$,
respectively.
The de\-Rham class $[\widetilde{\tau}]$ is the one corresponding to $[\tau]_\pi$
under the natural isomorphism between de Rham and Poisson cohomologies (see
(\ref{PoissonCohEq})). We can state 
\begin{lemma}
Let $M$ be a symplectic manifold, and let $\star, \star'$ be star products
with $C_1=C_1'=\frac{1}{2}\{\cdot,\cdot\}$. Then $[\widetilde{\tau}]
\in H_{\scriptscriptstyle {dR}}^2(M)$ depends only on $[\star],[\star']$.
\end{lemma}

Suppose now that $\star$ and $\star'$ satisfy $\star' \in
\Phi_{\scriptscriptstyle \mathcal E}([\star])$, where $\mathcal E = \Gamma^\infty(L)$
for a line bundle $L \to M$.
Let $(\mathcal{E}[[\lambda]],\bullet,\bullet')$ be a
bimodule quantization of $L$ corresponding to $\star, \star'$.
Let $R =R_1-R_1'$ be the  contravariant
connection on $L$ defined by $\bullet, \bullet'$,
and let $\Theta_R$ denote its curvature (see Appendix \ref{AppA}).

\begin{theorem}\label{CurvThm}
For $\; f,g \in C^\infty(M), \; s \in \mathcal E$, we have
$\tau(f,g) s = \Theta_R(df,dg)s$.
\end{theorem}
\begin{proof}
From (\ref{RelEqI}) we get, in second order,
\begin{equation} \label{I'}
R_2'(fg,s) + R_1'(C_1'(f,g),s) + C_2'(f,g)s = R_2'(f,gs) + R_1'(f,R_1'(g,s))
+ fR_2'(g,s).
\end{equation}
Similarly, from (\ref{RelEqII}) we get
\begin{equation} \label{II'}
R_2(s,fg) + R_1(s, C_1(f,g)) + s C_2(f,g) = R_2(sf,g) + R_1(R_1(s,f),g) +
R_2(s,f)g.
\end{equation}
Finally, from (\ref{RelEqIII}) we have
\begin{equation}\label{III'}
R_2(fs,g) + R_1(R_1'(f,s),g) + R_2'(f,s)g = R_2'(f,sg) + R_1'(f,R_1(s,g))
+ fR_2(s,g).
\end{equation}
Since we assume that $C_1=C_1'$,
subtracting (\ref{I'}) from (\ref{II'}) yields
\begin{eqnarray*}
 R(s,C_1(f,g)) -  R( R(s,f),g) +  R( R(s,g),f) & & \\
\mbox{} + (C_2-C_2')(f,g)s & = &
R_1'(g,R_1(s,f)) - R_1'(R_1'(g,s),f)\\
& & \mbox{} + R_2(sf,g) - R_2'(f,gs) \\
& & \mbox{} + R_2(s,f)g +  R(R_1(s,g), f)\\
& & \mbox{} + R(R_1'(s,f),g)- fR_2'(g,s) \\
& & \mbox{} + R_2'(fg,s) - R_2(s,fg).
\end{eqnarray*}
Using (\ref{III'}), we then get
\begin{eqnarray*}
 R(s,C_1(f,g)) -  R( R(s,f),g) +  R( R(s,g),f)& &\\
\mbox{}  + (C_2-C_2')(f,g)s & = &
R(R_1'(f,s),g)+ R(R_1(s,g),f)\\
& & \mbox{} + R_2'(fg,s) - R_2(s,fg) \\
& & \mbox{} + R_2(sf,g) + R_2(gs,f) \\
& & \mbox{} - R_2'(f,gs) - R_2'(g,sf).\\
\end{eqnarray*}
Taking the skew-symmetric part of this equation, and recalling that
$ \{f,g\} = C_1(f,g) - C_1(g,f)$, we finally have
$$
\tau(f,g) s = R(s,\{f,g\}) - R(R(s,f),g) + R(R(s,g),f).
$$
\end{proof}

Consider the natural map
\begin{equation}\label{NaturalEq}
i:H^2(M,\mathbb Z) \longrightarrow H^2_{\scriptscriptstyle{dR}}(M).
\end{equation}
We denote $H^2_{\scriptscriptstyle{dR}}(M,\mathbb Z):= i(H^2(M,\mathbb Z))$.
\begin{corollary}\label{ObsCor}
Suppose $\star$ and $\star'$ satisfy
$[\star'] \in \Phi_{\scriptscriptstyle \mathcal{E}}([\star])$, 
$\mathcal E = \Gamma^\infty(L)$.
Then $\frac{\im}{2\pi}[\tau]_{\pi}= \pch(L) \in
H^2_{\pi}(M,\mathbb Z)$, where $\pch(L)=\pi^*\cL(L)$ is the
Poisson-Chern class of $L$.
In particular, if $M$
is symplectic,
$\frac{\im}{2\pi}[\tilde{\tau}]=\cL(L) \in H^2_{\scriptscriptstyle{dR}}(M,\mathbb Z)$.
\end{corollary}

Corollary \ref{ObsCor} provides an integrality
 obstruction for Morita equivalent
star products on Poisson manifolds.
In the next three subsections, we will interpret these results in terms
of the characteristic classes of  star products.
\subsection{The symplectic case}\label{SympSub}

If $(M,\omega)$ is a symplectic manifold, the set of equivalence classes of star
products on $M$ can be described in terms of the second de Rham cohomology of
$M$ \cite{BCG97,Deligne95,Fed94a,NT95b,WX98}: There is a bijection
\begin{equation}\label{IdentSymEq}
\ch: \Def(M,\omega) \longrightarrow
[\omega] +  \lambda H^2_{\scriptscriptstyle {dR}}(M)[[\lambda]].
\end{equation}
\noindent The class $\ch(\star)$ is called the \emph{characteristic class} of $\star$.

The following result was proven
in \cite[Prop.~6.2]{GR99}.
\begin{lemma}\label{GuttLem}
Let $\star,\star'$ be star products on $M$, and let
$\widetilde{\tau}$ be the closed $2$-form defined in (\ref{TauTilEq}). Then
$[\widetilde{\tau}] = \frac{1}{\lambda}(\ch(\star')-\ch(\star))$ mod $\lambda$.
\end{lemma}

In order to study the semiclassical limit of
$\Phi:\Pic(M)\times\Def(M,\omega) \longrightarrow \Def(M,\omega)$,
let $S:H^2_{\scriptscriptstyle {dR}}(M)[[\lambda]] 
\longrightarrow H^2_{\scriptscriptstyle {dR}}(M)$ be the semiclassical
limit map $S(\sum_{r=0}^\infty[\omega_r]\lambda^r) = [\omega_1]$. With the
identification (\ref{IdentSymEq}), we may consider
$S:\Def(M,\omega) \longrightarrow H^2_{\scriptscriptstyle {dR}}(M)$.

Let $L \to M$ be a complex line bundle, and let $\mathcal E =
\Gamma^\infty(L)$.

\begin{theorem}\label{SemiClassThm}
The following diagram commutes:

\[
\begin{CD}
\Def(M,\omega) @>\Phi_{\scriptscriptstyle \mathcal E}>>
\Def(M,\omega)\\
@V{S}VV  @VV{S}V\\
H^2_{\scriptscriptstyle {dR}}(M) 
@>{\widehat{\Phi}_{\scriptscriptstyle \mathcal E}}>> H^2_{\scriptscriptstyle {dR}}(M),
\end{CD}
\]
where $\widehat{\Phi}_{\scriptscriptstyle \mathcal E}([\alpha])=
[\alpha] + \frac{2\pi}{\im}\cL(L)$.
\end{theorem}

\begin{proof}
The proof follows directly from Lemma \ref{GuttLem} and Corollary \ref{ObsCor}.
\end{proof}

Recall that $\Pic(M) \cong H^2(M,\mathbb Z)$ and
that the kernel of $i$ (see (\ref{NaturalEq}))
is given by the torsion elements in $H^2(M,\mathbb Z)$. We then have the
following

\begin{corollary}
Let $(M,\omega)$ be a symplectic manifold, and suppose $H^2(M,\mathbb Z)$
is free. Then
the action $\Phi: H^2(M,\mathbb Z)\times
\Def(M,\omega)\longrightarrow
\Def(M,\omega)$ is faithful. 
\end{corollary}

Recall that for star products $\star,\star'$ on $M$, their relative Deligne
class is defined by $t(\star,\star')=\ch(\star)-\ch(\star') \in \lambda
H^2_{\scriptscriptstyle {dR}}(M)[[\lambda]]$. 
We write $t(\star,\star')= \lambda t_0(\star,\star')
+ O(\lambda^2)$. We have the following immediate consequence of
Theorem \ref{SemiClassThm} phrased in terms of relative classes.

\begin{corollary}\label{ObsSymCor}
If $\star,\star'$ are Morita equivalent star products on a symplectic
manifold $(M,\omega)$, then there exists a symplectomorphism
$\psi:M\longrightarrow M$ such that $\frac{\im}{2\pi}t_0(\star,\psi^*(\star')) 
\in H^2_{\scriptscriptstyle {dR}}(M,\mathbb Z)$.
Conversely, for any star product $\star$
on $M$ and $[\alpha ]\in H^2_{\scriptscriptstyle {dR}}(M,\mathbb Z)$,
there is a star product $\star'$ Morita equivalent to $\star$ such that
$t(\star,\star')=\frac{2\pi}{\im}[\alpha]\lambda + O(\lambda^2)$.
\end{corollary}

\subsection{The Poisson case}\label{PoissSub}\label{PoissonSub}

For an arbitrary Poisson manifold $(M,\pi)$, Kontsevich
constructed in \cite{Kon97b} a bijection
\begin{equation}\label{PoissCarEq}
\ch: \Def(M,\pi) \longrightarrow
\{\pi_\lambda = \pi + \lambda \pi_1 + \ldots \in \chi^2(M)[[\lambda]],
[\pi_\lambda ,\pi_\lambda ]=0 \}/F,
\end{equation}
where $F$ is the group
$\{\exp(\sum_{r=1}^\infty D_r\lambda^r), \; D_r \in \Der(C^\infty(M))\}$, acting on
formal Poisson structures by $\qT(\pi_{\lambda})=\pi_{\lambda}'$ if and only if
$\pi_{\lambda}'(df,dg)= \qT^{-1}\pi_{\lambda}(d(\qT(f)),d(\qT(g)))$, for
$\qT \in F$. This correspondence is a result
of a more general fact \cite{Kon97b}:
there exists an $L_\infty$-quasi-isomorphism $\mathcal U$
from the graded Lie algebra of multivectors fields on $M$ (with zero
differential and Schouten bracket), $\mathfrak{g}_1$,
into the graded Lie algebra of
multidifferential operators on $M$ (with Hochschild differential and
Gerstenhaber bracket), $\mathfrak{g}_2$.
Given such an $\mathcal U$, for every
formal Poisson structure $\pi_\lambda$ we can define a star product
$\star_{\pi_\lambda }$ by
\begin{equation}\label{KonStarEq}
f \star_{\pi_\lambda } g := fg + \sum_{r=1}^\infty \frac{\lambda^r}{r!}
\mathcal{U}_r(\underbrace{\pi_{\lambda}\wedge\ldots \wedge \pi_{\lambda}}_r)
(f \otimes g),
\end{equation}
where $\mathcal{U}_r:\bigwedge^r \mathfrak{g}_1 \longrightarrow
\mathfrak{g}_2$ are the Taylor coefficients of $\mathcal U$. Moreover,
Kontsevich showed that one can choose
$\mathcal{U}_1: \mathfrak{g}_1 \longrightarrow \mathfrak{g}_2$
to be  the natural embedding of multivector fields into
multidifferential operators (note that this embedding does \emph{not}
preserve brackets).

If $\pi_{\lambda} =\pi + \lambda \pi_1 + \ldots$ is a formal
Poisson structure on $M$, the integrability equation $[\pi_{\lambda},
\pi_{\lambda}]=0$ immediately implies that $d_{\pi}\pi_1 = 0$.
\begin{lemma}\label{SemiPoissLem}
If $\pi_{\lambda} =\pi + \lambda \pi_1 + \ldots$ and
$\pi_{\lambda}' =\pi + \lambda \pi_1' + \ldots$ are equivalent formal Poisson
structures, then $[\pi_1]_{\pi} = [\pi_1']_{\pi}$.
\end{lemma}
\begin{proof}
Let $\qT=\exp(\sum_{r=1}^\infty D_r\lambda^r) \in F$.
A simple computation shows that if
$\qT(\pi_{\lambda})=\pi_{\lambda}'$, then
$$
\pi_1' = \pi_1 - d_{\pi}X_1,
$$
where $X_1 \in \chi(M)$ is defined by $\mathcal{L}_{X_1}=D_1$. Thus
$[\pi_1]_{\pi} = [\pi_1']_{\pi}$.
\end{proof}

With the identification given in (\ref{PoissCarEq}), we define the semiclassical
limit map 
$$
S: \Def(M,\pi) \longrightarrow H^2_\pi(M), \;\;
S([\pi_{\lambda}]) = [\pi_1]_\pi,
$$
where $[\pi_{\lambda}]$ is the equivalence class of the formal Poisson structure
$\pi_{\lambda} =\pi + \lambda \pi_1 + \ldots$.

\begin{lemma}\label{IndepLem}
Let $\star$ and $\star'$ be star products on $(M,\pi)$, with
$\ch(\star)=[\pi + \lambda \pi_1 + \ldots]$ and
$\ch(\star')=[\pi + \lambda \pi_1' + \ldots]$. Let
$\tau$ be as in Lemma \ref{TauLem}.
Then $[\tau]_\pi = [\pi_1]_{\pi} - [\pi'_1]_{\pi}$.
\end{lemma}
\begin{proof}
Since in our convention  $C_1^{\mbox{\tiny{skew}}}=\frac{1}{2}\{\cdot,\cdot\}$,
we use  Kontsevich's explicit construction for the formal Poisson structure
$\frac{1}{2}\pi_\lambda$.
The expression of Kontsevich's star products in terms of the maps $\mathcal{U}_r$
is
\begin{eqnarray*}
f\star_{\pi_{\lambda}}g& =&fg + \lambda\mathcal{U}_1(\frac{1}{2}
\pi_{\lambda})(f\otimes g)
+ \frac{\lambda^2}{2}
 \mathcal{U}_2(\frac{1}{2}\pi_{\lambda}\wedge
 \frac{1}{2}\pi_{\lambda})(f\otimes g) +
\ldots\\
& = & fg + \frac{\lambda }{2} \pi(df,dg) + \lambda^2(\frac{1}{2}\pi_1(df,dg)+
\frac{1}{8}\mathcal{U}_2(\pi \wedge \pi)(f\otimes g)) +
\ldots
\end{eqnarray*}
\noindent Since $\star$ is equivalent to $\star_{\pi_{\lambda}}$,
and $\star'$ is equivalent to $\star_{\pi_{\lambda}'}$, by Lemma \ref{IndepLem}
it suffices to compute $\tau$ for $\star_{\pi_{\lambda}}$ and
$\star_{\pi_{\lambda}'}$.
It is clear from the expression just above that $\tau=\pi_1 -\pi_1'$.
\end{proof}

Let $L\to M$ be a complex line bundle, and $\mathcal E = \Gamma^\infty(L)$.
The following result follows from Lemma \ref{IndepLem} and Theorem \ref{CurvThm}.
\begin{theorem}\label{SemiClassPoissThm}
The following diagram commutes:
\[
\begin{CD}
\Def(M,\pi) @>\Phi_{\scriptscriptstyle \mathcal E}>>
\Def(M,\pi)\\
@V{S}VV  @VV{S}V\\
H^2_{\pi}(M) @>{\widehat{\Phi}_{\scriptscriptstyle \mathcal E}}>> H^2_{\pi}(M),
\end{CD}
\]
where $\widehat{\Phi}_{\scriptscriptstyle \mathcal E}([\alpha])=
[\alpha] - \frac{2\pi}{\im}\pch(L) = [\alpha] - \frac{2\pi}{\im} {\pi^*\cL(L)}$.
\end{theorem}

Hence, for a star product $\star$ on $(M,\pi)$, each element in
$H^2_{\pi}(M,\mathbb Z)=\pi^*H^2_{\scriptscriptstyle{dR}}(M,\mathbb Z)$
corresponds to a different equivalence class of star products
Morita equivalent to $\star$. 
Theorem \ref{SemiClassPoissThm} shows, in particular, that the
semiclassical limit of
$\Phi$ is trivial when $\pi$ induces the trivial map in cohomology;
the case $\pi = 0$ will
be discussed in Section \ref{ZeroSub}.

A bivector field $\pi_1$ on a Poisson manifold $(M,\pi)$ is called an
\emph{infinitesimal deformation} of $\pi$ if $d_{\pi}\pi_1=0$. 

\begin{corollary}\label{PoissDefCor}
Suppose $\pi_1$ is an infinitesimal deformation that extends to  a
formal Poisson structure $\pi_{\lambda}$. Then the same holds for
$\pi_1 + \alpha$ if $\frac{\im}{2\pi}[\alpha]_{\pi} \in H^2_{\pi}(M,\mathbb Z)$.
\end{corollary}
\subsection{Deformations of the zero Poisson structure}
\label{ZeroSub}
As mentioned in the Section \ref{PoissonSub}, 
Theorem \ref{SemiClassPoissThm} does not provide
much information about the orbits of star products corresponding to the null
Poisson structure. We will show that the picture, in this case, is analogous
to  Sections  \ref{SympSub}, \ref{PoissonSub}, but 
in higher orders of $\lambda$.

Let $(M,\pi)$ be a Poisson manifold, with $\pi = 0$. For simplicity, we will
identify equivalence classes of
star products on $M$ with their characteristic classes as in (\ref{PoissCarEq}).
In order to understand the action of $\Phi$ on $\Def(M,\pi)$,
consider the disjoint union
\begin{equation}\label{UnionEq}
\Def(M,\pi) = \bigcup_{m\geq 1} \Def^m(M,\pi) \cup [0],
\end{equation}
where $[0]$ denotes the equivalence class of the trivial formal
Poisson structure on $M$, and $\Def^m(M,\pi)$ is the set of equivalence classes
of formal Poisson structures of the form $\pi_{\lambda}^m= \lambda^m(\pi_m + \lambda
\pi_{m+1} + O(\lambda^2))$, $\pi_m \neq 0$, $m\geq 1$.
Note that we can decompose each $\Def^m(M,\pi)$ into a disjoint union of sets
$\Def^m(M,\pi_m)$, given by equivalence classes of formal Poisson
structures of the form $\lambda^m(\pi_m + O(\lambda))$ for a \emph{fixed}
Poisson structure $\pi_m \neq 0$. We can always choose a star product
$\star = \sum_{r=0}^\infty C_r \lambda^r$ corresponding to a class in $\Def^m(M,\pi_m)$
with $C_1=C_2=\ldots=C_{m}=0$.
It is easy to check that all the results 
in the previous subsections of Section
\ref{EquivSec} hold for such star products, 
with a shift in order by $\lambda^m$. For instance,
the same arguments as in Theorem \ref{StarActThm} show that
$\Def^m(M,\pi_m)$ is invariant under $\Phi$.

\begin{corollary} The trivial class
$[0]$ is a fixed point for $\Phi$.
\end{corollary}

Let $S_m:\Def^m(M,\pi_m)\longrightarrow H^2_{\pi_m}(M)$ be defined by 
$S_m(\lambda^m(\pi_m + \lambda\pi_{m+1} + O(\lambda^2))) = [\pi_{m+1}]_{\pi_m}$. 
Let $L \to M$ be a line bundle and $\mathcal E = \Gamma^\infty(L)$.
Just as in Theorem \ref{SemiClassPoissThm}, one can show the following theorem. 
\begin{theorem}\label{ZeroSemiClassThm}
The following diagram commutes:
\[
\begin{CD}
\Def^m(M,\pi_m) @>\Phi_{\mathcal E}>> 
\Def^m(M,\pi_m)\\
@V{S_m}VV  @VV{S_m}V\\
H^2_{\pi_m}(M) @>{\widehat{\Phi}_{\mathcal E}}>> H^2_{\pi_m}(M),
\end{CD}
\]
where $\widehat{\Phi}_{\mathcal E}([\alpha])= 
[\alpha] - \frac{2\pi}{\im}(\pi_m^*\cL(L))$.
\end{theorem}

Hence, for a star product in $\Def^m(M,\pi_m)$, each element in 
$H^2_{\pi_m}(M,\mathbb Z)$
corresponds to an equivalence class of star products Morita equivalent to it.

\appendix
\section{Poisson cohomology, contravariant connections and Poisson-Chern classes}
\label{AppA}
Let $(M,\pi)$ be a Poisson manifold. The Poisson tensor $\pi$ defines a bundle
morphism
$$
\tilde{\pi}: T^*M \longrightarrow TM,\;\; \alpha \mapsto \pi(\cdot,\alpha),
$$
inducing a map on sections $\tilde{\pi}: \Omega^1(M)
\longrightarrow \chi(M)$. The vector field $\tilde{\pi}(df)=X_f$ is called
the \emph{hamiltonian vector field} of $f$.
We can use $\tilde{\pi}$ to define a Lie algebra bracket on $\Omega^1(M)$:
\begin{equation}\label{Lie1FormEq}
[\alpha,\beta] = -\mathcal{L}_{\tilde{\pi}(\alpha)}\beta +
\mathcal{L}_{\tilde{\pi}(\beta)}\alpha
- d(\pi(\alpha,\beta)), \;\; \alpha,\beta \in \Omega^1(M).
\end{equation}
\noindent The map  $-\tilde{\pi}:\Omega^1(M)
\longrightarrow \chi(M)$ is a Lie algebra homomorphism, and this
makes $T^*M$ into a \emph{Lie algebroid} (see \cite[Chp.~16]{SilWein99}).

The Poisson
tensor $\pi \in \chi^2(M)$ can be used to define a differential
\begin{equation}\label{PDiffEq}
d_{\pi}: \chi^k(M) \longrightarrow \chi^{k+1}(M),\; \;\; d_{\pi} = [\pi,\cdot],
\end{equation}
where $[\cdot,\cdot]$ is the Schouten bracket \cite{Vais94}.
\begin{definition}
The cohomology groups of the complex $(\chi^{\bullet},d_{\pi})$ are called the
{\it Poisson cohomology} groups of $M$ and denoted
$H^k_\pi(M)$.
\end{definition}
The map $\tilde{\pi}$ induces a map $\pi^*:\Omega^{\bullet(M)}
\longrightarrow \chi^{\bullet}(M)$  intertwining differentials, and
therefore gives rise to a morphism in cohomology
\begin{equation}\label{PoissonCohEq}
\pi^*: H^k_{\scriptscriptstyle {dR}}(M) \longrightarrow H^k_\pi(M),
\end{equation}
which is an isomorphism when $\pi$ is symplectic.
 We define integral (resp. real) Poisson cohomology
 as the image of integral (resp. real) de\-Rham cohomology classes on $M$
under $\pi^*$, i.e.,
$H^k_\pi(M,\mathbb Z)= \pi^*H^k_{\scriptscriptstyle {dR}}(M,\mathbb Z)$
(resp. $H^k_\pi(M,\mathbb R)= \pi^*H^k_{\scriptscriptstyle {dR}}(M,\mathbb R)$).

The key ingredient in defining contravariant connections on vector bundles over
Poisson manifolds
is to think of $T^*M$ as a ``new'' tangent bundle to $M$,
using its Lie algebroid structure
(see \cite{Rui2000}).

Let $E\to M$ be a complex vector bundle over a Poisson manifold $(M,\pi)$.
\begin{definition}
A \emph{contravariant connection}
on $E$
is a $\mathbb C$-linear map 
$D: \Gamma^\infty(E)\otimes \Omega^1(M) \longrightarrow \Gamma^\infty(E)$
so that
\begin{enumerate}
\item $D_{f \alpha}s = f D_\alpha s$
\item $D_\alpha(f s)= f D_\alpha s + \alpha(X_f) s$,
\end{enumerate}
for $\alpha \in \Omega^1(M)$,
$f \in C^\infty(M)$.
The \emph{curvature} of a contravariant connection $D$ is the map
$\Theta_D:\Omega^1(M)\otimes\Omega^1(M) \longrightarrow \End(\Gamma^\infty(E))$,
$$
\Theta_D(\alpha,\beta) s = D_\alpha D_\beta s - D_\beta D_\alpha s +
D_{[\alpha, \beta]}s.
$$
\end{definition}

It is easy to see that, if $\nabla$ is any connection (in the usual sense) on
$E$, then it induces a contravariant connection by $D_{df}=\nabla_{X_f}$.
On symplectic manifolds this is the only way that contravariant connections can arise.
Thus this notion is mostly important in degenerate situations.

A bilinear map
$D':\Gamma^\infty(E) \times C^\infty(M) \longrightarrow \Gamma^\infty(E)$,
satisfying
\begin{eqnarray}
D'(s, f\cdot g) & = & D'(s,f)g + D'(s,g)f, \label{ContConnEqI}\\
D'(s\cdot f,g) & = & D'(s,g)f + s \{f,g\},\label{ContConnEqII}
\end{eqnarray}
provides an equivalent definition of a contravariant connection.
The definitions are related by the formula
$$
D'(s,f)=D_{df}s.
$$

If $E=L \to M$ is a line bundle, then the curvature $\Theta_D$
of a contravariant
connection  defines a bivector field on $M$, closed
with respect to $d_\pi$ \cite{Vais91}. As in the case of usual connections,
its Poisson cohomology is a well-defined class, independent of the connection.
\begin{definition}
let $D$ be a contravariant connection on a line bundle $L \to M$, and let
$\Theta_D$ be its curvature.
We call the class
$\pch(L):=\frac{\im}{2\pi}[\Theta_D]_\pi \in H^2_\pi(M)$
the \emph{Poisson-Chern class} of $L$.
\end{definition}
\noindent It is clear that $\pch(L)=\pi^*(\cL(L))$.

\begin{footnotesize}

\end{footnotesize}


\begin{thebibliography}{10}

\bibitem {Bass68}
{\sc Bass, H.: }\newblock {\em Algebraic ${K}$-theory}.
\newblock W. A. Benjamin, Inc., New York-Amsterdam, 1968.

\bibitem {BFFLS78}
{\sc Bayen, F., Flato, M., Fr{{\o}}nsdal, C., Lichnerowicz, A., Sternheimer,
  D.: }\newblock {\em Deformation Theory and Quantization}.
\newblock Ann. Phys.  {\bf 111} (1978), 61--151.


\bibitem {BCG97}
{\sc Bertelson, M., Cahen, M., Gutt, S.: }\newblock {\em Equivalence of Star
  Products}.
\newblock Class. Quantum Grav.  {\bf 14} (1997), A93--A107.

\bibitem {Bu2001b}
{\sc Bursztyn, H.: }\newblock {\em Poisson vector bundles, contravariant connections and deformations}.
\newblock Contribution to Proceedings of
the Workshop on Noncommutative Geometry and String
Theory, Keio University, Japan, March 2001.


\bibitem {BuWa2001}
{\sc Bursztyn, H., Waldmann, S.: }\newblock {\em The characteristic classes 
of Morita equivalent star products on symplectic manifolds}.
\newblock Preprint math.QA/0106178.

\bibitem {BuWa2000}
{\sc Bursztyn, H., Waldmann, S.: }\newblock {\em {$^*$}-Ideals and Formal
  Morita Equivalence of {$^*$}-Algebras}.
\newblock Internat. J. Math. {\bf 12}.5 (2001), 555 -- 577. 
\newblock Preprint math.QA/0005227.

\bibitem {BuWa2000b}
{\sc Bursztyn, H., Waldmann, S.: }\newblock {\em Deformation Quantization of
  Hermitian Vector Bundles}.
\newblock Lett. Math. Phys.  {\bf 53}.4 (2000), 349--365.
\newblock Preprint math.QA/0009170.

\bibitem {BuWa99a}
{\sc Bursztyn, H., Waldmann, S.: }\newblock {\em Algebraic Rieffel Induction,
  Formal Morita Equivalence and Applications to Deformation Quantization}.
\newblock J. Geom. Phys.  {\bf 37}.4 (2001), 307--364.
\newblock Preprint math.QA/9912182.

\bibitem {SilWein99}
{\sc Cannas~da Silva, A., Weinstein, A.: }\newblock {\em Geometric models for
  noncommutative algebras}.
\newblock American Mathematical Society, Providence, RI, 1999.

\bibitem {ConnDoSch98}
{\sc Connes, A., Douglas, M.~R., Schwarz, A.: }\newblock {\em Noncommutative
  geometry and matrix theory: compactification on tori}.
\newblock J. High Energy Phys.  {\bf 1998}.2, Paper 3, 35 pp. (electronic).

\bibitem {Deligne95}
{\sc Deligne, P.: }\newblock {\em D\'eformations de l'alg\`ebre des fonctions
  d'une vari\'et\'e symplectique: comparaison entre {F}edosov et {D}e {W}ilde,
  {L}ecomte}.
\newblock Selecta Math. (N.S.)  {\bf 1}.4 (1995), 667--697.


\bibitem {DS2000}
{\sc Dito, G., Sternheimer, D. (eds.): }\newblock {\em Conf\'erence Moshe Flato
  1999: Quantization, Deformations, Symmetries}.
\newblock {\em Mathematical Physics Studies} no. {\bf 23}.
\newblock Kluwer Academic Press, Dordrecht, Boston, London, 2000.

\bibitem {EmmWe96}
{\sc Emmrich, C., Weinstein, A.: }\newblock {\em Geometry of the transport
  equation in multicomponent {W}{K}{B} approximations}.
\newblock Comm. Math. Phys.  {\bf 176}.3 (1996), 701--711.

\bibitem {Fed94a}
{\sc Fedosov, B.~V.: }\newblock {\em A Simple Geometrical Construction of
  Deformation Quantization}.
\newblock J. Diff. Geom.  {\bf 40} (1994), 213--238.

\bibitem {Fed96}
{\sc Fedosov, B.~V.: }\newblock {\em Deformation Quantization and Index
  Theory}.
\newblock Akademie Verlag, Berlin, 1996.

\bibitem {Rui2000}
{\sc Fernandes, R.~L.: }\newblock {\em Connections in {P}oisson geometry. {I}.
  {H}olonomy and invariants}.
\newblock J. Differential Geom.  {\bf 54}.2 (2000), 303--365.
\newblock math.DG/0001129.

\bibitem {Ger64}
{\sc Gerstenhaber, M.: }\newblock {\em On the Deformation of Rings and
  Algebras}.
\newblock Ann. Math.  {\bf 79} (1964), 59--103.

\bibitem {GS88}
{\sc Gerstenhaber, M., Schack, S.~D.: }\newblock {\em Algebraic Cohomology and
  Deformation Theory}.
\newblock In: {\sc Hazewinkel, M., Gerstenhaber, M. (eds.): }\newblock {\em
  Deformation Theory of Algebras and Structures and Applications},   13--264.
  Kluwer Academic Press, Dordrecht, 1988.

\bibitem {Gutt2000}
{\sc Gutt, S.: }\newblock {\em Variations on deformation quantization}.
\newblock In: {\sc Dito, G., Sternheimer, D. (eds.): }\newblock {\em
  Conf\'erence Moshe Flato 1999: Quantization, Deformations, Symmetries}.
  \cite{DS2000}.
\newblock  math.DG/0003107.

\bibitem {GR99}
{\sc Gutt, S., Rawnsley, J.: }\newblock {\em Equivalence of star products on a
  symplectic manifold; an introduction to Deligne's \v{C}ech cohomology
  classes}.
\newblock J. Geom. Phys.  {\bf 29} (1999), 347--392.

\bibitem {Hirz95}
{\sc Hirzebruch, F.: }\newblock {\em Topological methods in algebraic
  geometry}.
\newblock Springer-Verlag, Berlin, 1995.
\newblock Translated from the German and Appendix One by R. L. E.
  Schwarzenberger, With a preface to the third English edition by the author
  and Schwarzenberger, Appendix Two by A. Borel, Reprint of the 1978 edition.

\bibitem {JSW2001}
{\sc Jurco, B., Schupp, P., Wess, J.: }
\newblock {\em Noncommutative line bundle and Morita equivalence }.
\newblock Preprint hep-th/0106110.


\bibitem {Kon97b}
{\sc Kontsevich, M.: }\newblock {\em Deformation Quantization of Poisson
  Manifolds, I}.
\newblock Preprint  q-alg/9709040.

\bibitem {Lam99}
{\sc Lam, T.~Y.: }\newblock {\em Lectures on modules and rings}.
\newblock Springer-Verlag, New York, 1999.

\bibitem {Land2000}
{\sc Landsman, N.~P.: }\newblock {\em Quantized reduction as a tensor product}.
\newblock Preprint  math-ph/0008004.

\bibitem {Morita58}
{\sc Morita, K.: }\newblock {\em Duality for modules and its applications to
  the theory of rings with minimum condition}.
\newblock Sci. Rep. Tokyo Kyoiku Daigaku Sect. A  {\bf 6} (1958), 83--142.

\bibitem {NT95b}
{\sc Nest, R., Tsygan, B.: }\newblock {\em Algebraic Index Theorem for
  Families}.
\newblock Adv. Math.  {\bf 113} (1995), 151--205.

\bibitem {PoRe86}
{\sc Porta, H., Recht, L.: }\newblock {\em Classification of linear
  connections}.
\newblock J. Math. Anal. Appl.  {\bf 118}.2 (1986), 547--560.

\bibitem {RVW96}
{\sc Reshetikhin, N., Voronov, A., Weinstein, A.: }\newblock {\em Semiquantum
  geometry}.
\newblock J. Math. Sci.  {\bf 82}.1 (1996), 3255--3267.
\newblock Algebraic geometry, 5.


\bibitem {Rief89}
{\sc Rieffel, M.~A.: }\newblock {\em Deformation quantization of {H}eisenberg
  manifolds}.
\newblock Comm. Math. Phys.  {\bf 122}.4 (1989), 531--562.

\bibitem {RiefSch}
{\sc Rieffel, M.~A., Schwarz, A.: }\newblock {\em Morita equivalence of
  multidimensional noncommutative tori}.
\newblock Internat. J. Math.  {\bf 10}.2 (1999), 289--299.


\bibitem {Sch98}
{\sc Schwarz, A.: }\newblock {\em Morita equivalence and duality}.
\newblock Nuclear Phys. B  {\bf 534}.3 (1998), 720--738.

\bibitem {Ste98}
{\sc Sternheimer, D.: }\newblock {\em Deformation Quantization: Twenty Years
  After}.
\newblock Preprint  math.QA/9809056.

\bibitem {Vais91}
{\sc Vaisman, I.: }\newblock {\em On the geometric quantization of {P}oisson
  manifolds}.
\newblock J. Math. Phys.  {\bf 32}.12 (1991), 3339--3345.

\bibitem {Vais94}
{\sc Vaisman, I.: }\newblock {\em Lectures on the geometry of {P}oisson
  manifolds}.
\newblock Birkh\"auser Verlag, Basel, 1994.

\bibitem {Watts60}
{\sc Watts, C.~E.: }\newblock {\em Intrinsic characterizations of some additive
  functors}.
\newblock Proc. Amer. Math. Soc.  {\bf 11} (1960), 5--8.


\bibitem {Wei94}
{\sc Weinstein, A.: }\newblock {\em Deformation Quantization}.
\newblock S\'eminaire Bourbaki 46\`eme ann\'ee  {\bf 789} (1994).

\bibitem {WX98}
{\sc Weinstein, A., Xu, P.: }\newblock {\em Hochschild cohomology and
  characterisic classes for star-products}.
\newblock In: {\sc Khovanskij, A., Varchenko, A., Vassiliev, V. (eds.):
  }\newblock {\em Geometry of differential equations. Dedicated to V. I. Arnold
  on the occasion of his 60th birthday},   177--194. American Mathematical
  Society, Providence, 1998.

\end{thebibliography}
\end{document}